\documentclass[12pt]{amsart}
\usepackage{t1enc}
\usepackage[OT2,T1]{fontenc}
\usepackage{enumerate}

\newif\ifdeveloping

\newenvironment{cyr}{\fontencoding{OT2}\fontfamily{cmr}\selectfont}{}

\let\QED\qed
\newcommand{\prlabel}[1]{\renewcommand{\qed}{\QED${}_{\ref{#1}}$}}

\newtheorem{theorem}{Theorem}[section]
\newtheorem{lemma}[theorem]{Lemma}
\newtheorem{corollary}[theorem]{Corollary}
\newtheorem{claim}{Claim}[theorem]
\newtheorem{problem}[theorem]{Problem}

\newtheorem{qtheorem}{Theorem}

\theoremstyle{definition}
\newtheorem{definition}[theorem]{Definition}

\theoremstyle{remark}
  
{}

\ifdeveloping
\usepackage[notref,notcite]{showkeys}
\fi

\newcommand{\prtime}{{\count0=\time\divide\count0 by 60
\count1=-\count0\multiply\count1 by 60
\advance\count1 by \time
\the\count0:\the\count1}
}

\def\myheads#1;#2;{
\pagestyle{myheadings}
\markboth{{\sc\hfill #1\hfill\protect\makebox[0cm][r]{\rm\today; \prtime}}}
{{\sc\protect\makebox[0cm][l]{\rm\today;\ \prtime}\hfill #2\hfill}}
\thispagestyle{myheadings}
}
\newcommand{\acal}{{\mathcal A}}
\newcommand{\bcal}{{\mathcal B}}
\newcommand{\dcal}{{\mathcal D}}
\newcommand{\gcal}{{\mathcal G}}
\newcommand{\ical}{{\mathcal I}}
\newcommand{\pcal}{{\mathcal P}}
\newcommand{\ucal}{{\mathcal U}}
\newcommand{\vcal}{{\mathcal V}}
\newcommand{\ycal}{{\mathcal Y}}
\newcommand{\setm}{\setminus}
\newcommand{\empt}{\emptyset}
\newcommand{\subs}{\subset}
\def\<{\left\langle}
\def\>{\right\rangle}
\def\cf{\operatorname{cf}}
\def\br#1;#2;{\bigl[ {#1} \bigr]^ {#2} }
\newcommand{\ps}{\operatorname{ps}}
\newcommand{\ee}{\operatorname{e}}
\newcommand{\ext}{\operatorname{ext}}
\newcommand{\pp}{\operatorname{p}}
\newcommand{\pe}{\operatorname{pe}}

\newcommand{\tip}{\operatorname{tp}}
\newcommand{\rs}{\operatorname{rs}}
\theoremstyle{plain}

\def\cl#1{\overline{#1}}
\def\celh{\operatorname{\widehat c}}
\def\hh{\operatorname{H}}
\def\aa#1#2{A^{#1}_{#2}}
\def\aaa#1{A_{#1}}
\def\hlh{\operatorname{\widehat h}}
\def\ns{\operatorname{NS}}
\def\spre{\operatorname{s}}
\def\spreh{\operatorname{\widehat{s}}}
\begin{document}

\author[I. Juh\'asz]{Istv\'an Juh\'asz}
\address{Alfr{\'e}d R{\'e}nyi Institute of Mathematics}
\email{juhasz@renyi.hu}

\author[L. Soukup]{
Lajos Soukup }
\address{
Alfr{\'e}d R{\'e}nyi Institute of Mathematics }
\email{soukup@renyi.hu}

\author[Z. Szentmikl\'ossy]{
Zolt\'an Szentmikl\'ossy}
\address{E\"otv\"os University of Budapest}
\email{zoli@renyi.hu}

\subjclass[2000]{ }
\title{Resolvability of spaces having small spread or extent}
\thanks{The preparation of this paper was supported by OTKA grant no. 37758}

\begin{abstract}
In a recent paper O. Pavlov proved the following two interesting
resolvability results:
\begin{enumerate}
\item If a space $X$ satisfies $\Delta(X) > \ps(X)$ then
$X$ is maximally resolvable.
\item If a $T_3$-space $X$ satisfies $\Delta(X) > \pe(X)$ then
$X$ is $\omega$-resolvable.
\end{enumerate}
\smallskip Here $\ps(X)$ ($\pe(X)$) denotes the smallest successor
cardinal such that $X$ has no discrete (closed discrete) subset of
that size and $\Delta(X)$ is the smallest cardinality of a non-empty
open set in $X$. In this note we improve (1) by showing that
$\Delta(X)
>$ $\ps(X)$ can be relaxed to $\Delta(X) \ge$ $\ps(X)$. In particular,
if $X$ is a space of countable spread with $\Delta(X) > \omega$ then
$X$ is maximally resolvable.

The question if an analogous improvement of (2) is valid remains
open, but we present a proof of (2) that is simpler than Pavlov's.

\end{abstract}

\maketitle \ifdeveloping \myheads{Resolvability };{Resolvability };
\fi

\section{Introduction}

Given a cardinal ${\kappa}>1$, a topological space  is called {\em
${\kappa}$-resolvable} iff it contains ${\kappa}$ many disjoint
dense subsets. Denoting by ${\tau}^*(X)$ the family of nonempty open
subsets of a topological space $X$, we say that the space  $X$ is
{\em maximally resolvable } iff it is
$\operatorname{\Delta}(X)$-resolvable, where 
$\operatorname{\Delta}(X) = \min\bigl\{|G|:G\in {\tau}^*(X)\bigr\}$
is the so-called {\em dispersion character} of $X$. A space is
called $(<{\kappa})$-resolvable iff it is $\mu$-resolvable for all
$\mu < \kappa$. In this introduction we shall give three lemmas that
provide sufficient conditions for $\kappa$-resolvability. Finally, a
space that is not $\kappa$-resolvable is also called {\em
$\kappa$-irresolvable}.

El'kin   proved in \cite{El} that, for any cardinal ${\kappa}$,
every space may be written as the disjoint union of a hereditarily
${\kappa}$-irresolvable open subset and a
 ${\kappa}$-resolvable closed subset.
As Pavlov observed in the introduction of \cite{Pa}, this statement
 has the following reformulation.

\begin{lemma}\label{lm:elkin}
A topological space $X$ is ${\kappa}$-resolvable iff every nonempty
open subspace of $X$ includes a nonempty ${\kappa}$-resolvable
subset, in other words: iff $X$ has a $\pi$-network consisting of
${\kappa}$-resolvable subsets.
\end{lemma}

\newcommand{\ls}{\operatorname{ls}}
For any topological space $X$ we let $\ls(X)$ denote the minimum
number of left-separated subspaces needed to cover $X$. The
following lemma is implicit in the
 proof of \cite[Theorem 2.8]{Pa} and easily follows from the
 well-known fact that every space has a {\em dense} left-separated subspace,
 see e. g. \cite[2.9.c]{J}.

\begin{lemma}\label{lm:ls}
If for each $U\in {\tau}^*(X)$ we have $\ls(U)\ge {\kappa}$, that is
no nonempty open set in $X$ can be covered by fewer than $\kappa$
many left separated sets, then $X$ is ${\kappa}$-resolvable.
\end{lemma}

Our next lemma generalizes propositions 2.3 and 3.3 from \cite{Pa}.
We believe that our present approach is not only more general but
also simpler than that in \cite{Pa}. To formulate the lemma, we need
to introduce a piece of notation.

Given a family of sets $\acal$ and a cardinal ${\kappa}$, we denote
by $S_{\kappa}(\acal)$ the collection of all {\em disjoint}
subfamilies of $\acal$ of size less than $\kappa$, i. e.
\begin{displaymath}
S_{\kappa}(\acal)=\{\acal'\in \br \acal;<{\kappa}; : \acal' \mbox{
is disjoint}\}.
\end{displaymath}

\begin{lemma}\label{lm:jcal}
Let us be given a topological space $X$, a dense set $D \subs X$, an
infinite cardinal $\kappa \ge |D|$, moreover a family $\ical \subs
\pcal(X)$ of subsets of $X$. If for each $ x\in D$ and for any
$\ycal \in S_{\kappa}(\ical)$ there is a set $ Z\in \ical$  such
that $ \cup \ycal \cap Z = \emptyset$ and   $x\in \overline Z$ then
$X$ is ${\kappa}$-resolvable.
\end{lemma}

\begin{proof}
Let  $\{x_{\alpha}:{\alpha}<{\kappa}\} = D$ be a ${\kappa}$-abundant
enumeration of $D$, that is for any point $x \in D$ we have $a_x =
\{\alpha : x_\alpha = x\} \in [\kappa]^\kappa$. By a straightforward
transfinite recursion on ${\alpha}<{\kappa}$ we may then choose sets
$Z_{\alpha}\in \ical\cap \pcal(X\setm \cup_{{\nu}<{\alpha}}Z_{\nu})$
with  $x_{\alpha}\in \overline{Z_{\alpha}}$ for all $\alpha <
\kappa$. (Note that we have $\{Z_\nu : \nu < \alpha\} \in
S_\kappa(\ical)$ along the way.)

For any ordinal $i < \kappa$ and for any point $x \in D$ let
$\alpha^x_i$ be the $i$th element of the set $a_x$ and set
$$D_i=\bigcup\{Z_{\alpha^x_i}:x \in D\}.$$ Then clearly $D \subs
\overline{ D_i}$ , hence $\{D_i:i<{\kappa}\}$ is a disjoint family
of dense sets, witnessing that $X$ is ${\kappa}$-resolvable.
\end{proof}

As an illustration, note that if $|X|=\Delta(X)={\kappa}>{\lambda}$
and $t(x,X) \le \lambda$ holds for all points $x \in D$ of a set $D$
which is dense in the space $X$, then $D$, $ X$, $\kappa$, and
$\ical=\br X;\le {\lambda};$ satisfy the conditions of lemma
\ref{lm:jcal} and so $X$ is ${\kappa}$-resolvable. Thus we obtain
the following result as an immediate corollary of lemma
\ref{lm:jcal}.

\begin{corollary}
If $\Delta(X)> \sup\{ t(x,X) : x \in D\}$ for some dense set $D
\subs X$ then $X$ is maximally resolvable. In particular, if
$\Delta(X) > t(X)$ then $X$ is maximally resolvable.
\end{corollary}
The second statement is a theorem of Pytkeev from \cite{Py}.

\section{Improving Pavlov's result concerning spread}

As was mentioned in the abstract, in \cite{Pa} Pavlov defined
$\ps(X)$ as the smallest successor cardinal such that $X$ has no
discrete subset of that size. We recall from \cite[1.22]{J} the
related definition of $\widehat{s}(X)$ that is the smallest
uncountable cardinal such that $X$ has no discrete subset of that
size. Clearly, one has $\widehat{s}(X) \le$ $\ps(X)$ and
$\widehat{s}(X) =$ $\ps(X)$ iff $\widehat{s}(X)$ is a successor.
Finally, let us define $\operatorname{rs}(X)$ as the smallest uncountable regular
cardinal such that $X$ has no discrete subset of that size. Then we
have $\widehat{s}(X) \le\operatorname{rs}(X) \le$ $\ps(X)$ and 
$\widehat{s}(X) =\operatorname{rs}(X)$ iff $\widehat{s}(X)$ is regular.

In \cite{Pa} it was shown that if a space $X$ satisfies $\Delta(X)
>$ $\ps(X)$ then $X$ is maximally (i. e. $\Delta(X)$) resolvable. The
aim of this section is to improve this result by showing that the
assumption $\Delta(X) >\ps(X)$ can be relaxed to $\Delta(X) \ge $
$\operatorname{rs}(X)$.

Before doing that, however, we have to give an auxiliary result that
involves the cardinal function $\operatorname{h}(X)$, or more precisely its
"hatted" version $\hlh(X)$. We recall that ${\hlh}(X)$ is the
smallest uncountable cardinal such that $X$ has no right separated
subset of that size, or equivalently, the smallest uncountable
cardinal $\kappa$ with the property that any family $\ucal$ of open
sets in $X$ has a subfamily $\vcal$ of size $< \kappa$ such that
$\cup \vcal = \cup \ucal$, see e. g. \cite[2.9.b]{J}.

\begin{lemma}\label{lm:h}
 If ${\kappa}$ is an uncountable regular cardinal and
$$|X|\ge {\kappa}\ge \hlh(X)$$ then $X$ contains a
${\kappa}$-resolvable subspace $X^*$.
\end{lemma}

\begin{proof}
We can assume without loss of generality that
$X=\<{\kappa},{\tau}\>$. Let us denote by $\ns(\kappa)$ the ideal of
non-stationary subsets of $\kappa$ and set
$\gcal=\{U\in{\tau}:U\in\ns({\kappa})\}$. Since $\hlh(X)\le
{\kappa}$ there is $\gcal'\in \br\gcal;<{\kappa};$ with
$\cup\gcal'=\cup\gcal = G$. Then $G \in \ns({\kappa})$ because the
ideal $\ns({\kappa})$ is ${\kappa}$-complete.

Let us now consider the set
\begin{displaymath}
T=\{x\in {\kappa}:\exists C_x\subs {\kappa}\text{ club }
(\forall S\subs C_x \text{ if }S\in \ns({\kappa})
\text{ then } x\notin \overline S)\}.
\end{displaymath}
\begin{claim}
$T\in \ns({\kappa})$.
\end{claim}
Assume, on the contrary, that $T$ is stationary in $\kappa$. Fix for
each $x \in T$ a club $C_x$ as above. Then the diagonal intersection
$$C=\bigtriangleup\{C_x : x \in T\}$$ is again club and so
 $C\cap T$ is stationary in ${\kappa}$ as well. We may then choose a
set $S\in \br C\cap T;{\kappa};$ that is non-stationary. But then
for each $x\in S$ we have $$S\setm (x+1)\subs C\setm  (x+1)\subs
C_x,$$ hence by the choice of $C_x$ we have $x\notin
\overline{S\setm (x+1)}$. Consequently, $S$ is right separated in
its natural well-ordering, contradicting the assumption $\hlh(X)\le
{\kappa}$, and so our claim has been verified.

Finally, put $X^*=X\setm (G\cup T)$ and  $\ical=\ns({\kappa})\cap
\pcal(X^*)$. Then lemma \ref{lm:jcal} can be applied to the space
$X^*$, with itself as a dense subspace, the cardinal $\kappa$, and
the family $\ical$. Indeed, for any point $x \in X^*$ and for any
non-stationary set $Y \subs X^*$ there is a club set $C \subs X^*
\backslash Y$, and then $x \notin T$ implies that $x \in
\overline{Z}$ for some non-stationary set $Z \subs C$. (We have, of
course, used here that $\ical$ is $\kappa$-complete.) This shows
that $X^*$ is indeed ${\kappa}$-resolvable.
\end{proof}

We are now ready to formulate and prove the promised improvement of
Pavlov's theorem.

\begin{theorem}\label{tm:spread}
Let $X$ be a space and $\kappa$ be a regular cardinal such that
$$\widehat{s}(X) \le \kappa \le \Delta(X),$$ then $X$ is $\kappa$-resolvable.
Consequently, if $\Delta(X) \ge \rs(X)$ holds for a space $X$ then
$X$ is maximally resolvable. In particular, any space of countable
spread and uncountable dispersion character is maximally resolvable.
\end{theorem}

\begin{proof}

In view of lemma \ref{lm:elkin} it suffices to show that any
non-empty open subset $G$ of $X$ includes a $\kappa$-resolvable
subspace. To this end, note that, trivially, for each $G\in
{\tau}^*(X)$ we have either
  \begin{enumerate}[(i)]
   \item  $\ls(H)\ge {\kappa}$ for all $H\in{\tau}^*(G)$,

  \end{enumerate}
or
\begin{enumerate}[(i)]\addtocounter{enumi}{1}
\item  $\ls(H)<
  {\kappa}$ for some $H\in{\tau}^*(G)$.
\end{enumerate}
In case (i) $G$ itself is ${\kappa}$-resolvable by lemma
\ref{lm:ls}. In case (ii)  we claim that $\hlh(H)\le {\kappa}$ holds
true and therefore $H$ (and hence $G$) contains a
${\kappa}$-resolvable subset by lemma \ref{lm:h}. Assume, on the
contrary, that $R\subs H$ is right-separated and has cardinality
${\kappa}$. Since $H = \bigcup\{L_{\alpha}:{\alpha}<\ls(H)\}$, where
the sets $L_{\alpha}$ are all left-separated,  there is an
${\alpha}<\ls(H) < \kappa$ such that $|R\cap L_{\alpha}| = {\kappa}$
because ${\kappa}$ is regular. But then the subspace $R\cap
L_{\alpha}$ is both right and left separated, hence (see e. g.
\cite[2.12]{J})  it contains a discrete subset of size $|R\cap
L_{\alpha}|={\kappa}$, contradicting our assumption that
$\spreh(X)\le {\kappa}$.

If $\Delta(X)$ is regular then this immediately yields that $X$ is
maximally resolvable, while if $\Delta(X)$ is singular then, as
$\rs(X)$ is regular, we have
$$\Delta(X)
> \rs(X)^+ \ge \ps(X),$$
hence Pavlov's result \cite[2.9]{Pa} may be applied to get the
second part, of which the third is a special case.

\end{proof}

It is natural to raise the question if theorem \ref{tm:spread} could
be further improved by replacing $\rs(X)$ with $\widehat{s}(X)$ in
it. Of course, this is really a problem only in the case when
$$\Delta(X) = \widehat{s}(X) = \lambda$$ is a singular cardinal.
Recall now that Hajnal and Juh\'asz proved in \cite{HJ} (see also
\cite[4.2]{J}) that $\spreh(X)$ can not be singular {\em strong
limit} for a Hausdorff space $X$. Consequently, the above mentioned
strengthening is valid for Hausdorff spaces provided that all
singular cardinals are strong limit, in particular if GCH holds.

\begin{corollary}
Assume that for every (infinite) cardinal $\kappa$ the power
$2^\kappa$ is a finite successor of $\kappa$ (or equivalently, all
singular cardinals are strong limit). Then every Hausdorff space $X$
satisfying $\Delta(X)\ge\spreh(X)$ is maximally resolvable.
\end{corollary}

It is also known (see e. g. \cite[4.3]{J}) that $\spreh(X)$ can not
have countable cofinality for a strongly Hausdorff, in particular
for a $T_3$ space $X$. Hence the first interesting ZFC question that
is left open by theorem \ref{tm:spread} is the following.

\begin{problem}

Assume that  $X$ is a $T_3$ space satisfying
$$\spreh(X)=\operatorname{\Delta}(X) = \aleph_{\omega_1} .$$ Is $X$
then (maximally) resolvable?
\end{problem}

 It is clear that if in
theorem \ref{tm:spread} we have $\Delta(X) = \lambda > \rs(X)$ then
the first part may be applied to any regular cardinal $\kappa$ with
$\rs(X) \le \kappa \le \lambda$, hence if $\lambda$ is singular then
we obtain that $X$ is $(<\lambda)$-resolvable without any reference
to Pavlov's result. This is of significance because the proof of
Pavlov's theorem in the case when $\Delta(X)$ is singular is rather
involved. However, if in addition $\lambda$ has countable cofinality
then no reference to Pavlov's proof is needed because of the
following result of Bhaskara Rao.

\begin{qtheorem}[Bhaskara Rao, \cite{Ba}]
If $\cf({\lambda})={\omega}$ and the space $X$ is
  $(<\lambda)$-resolvable  then
$X$ is also ${\lambda}$-resolvable.
\end{qtheorem}

The question if the analogous result can be proved for singular
cardinals of uncountable cofinality is one of the outstanding open
problems in the area of resolvability and was already formulated in
\cite{df}. We just repeat it here.

\begin{problem}\label{pr:br}
Assume that  ${\lambda}$ is a singular cardinal with
$\cf({\lambda})>{\omega}$ and the space $X$ is
$(<\lambda)$-resolvable. Is it true then that $X$ is also
${\lambda}$-resolvable?
\end{problem}

We close this section by giving a partial affirmative answer to
problem \ref{pr:br}. At the same time we shall also show how the
first part of theorem \ref{tm:spread} implies the second in case
$\Delta(X)$ is singular, thus making our proof of \ref{tm:spread}
self-contained. To do this, we shall first fix some notation.

\begin{definition}
For any space $X$ we let $\dcal(X)$ denote the family of all dense
subsets of $X$. Next, we set $$\mathcal{F}(X) = \cup\{\mathcal{D}(U)
: U \in \tau^*(X)\};$$ we call the members of $\mathcal{F}(X)$, i.
e. dense subsets of (non-empty) open sets, {\em fat} sets in $X$.

For a subspace $Y\subs X$ and  a cardinal $\nu$ we let
$$\mathcal{H}(Y,\nu) = \mathcal{F}(X) \cap [Y]^{\le \nu},$$
in other words, $\mathcal{H}(Y,\nu)$ is the family of all fat (in
$X$ !) subsets of $Y$ of size at most $\nu$. It is easy to see that
if $c(X) \le \nu$ and $\mathcal{H}(Y,\nu)$ is non-empty then there
is a member $H(Y,\nu) \in \mathcal{H}(Y,\nu)$ of maximal closure, i.
e. such that
$$\overline{H(Y,\nu)} = \overline{\cup \mathcal{H}(Y,\nu)}.$$ (If
$\mathcal{H}(Y,\nu)$ is empty then we set $H(Y, \nu) = \emptyset$.)
Clearly, if $Y \subs Z \subs X$ and $c(X) \le \nu$ then we have
$$\overline{H(Y,\nu)} \subs \overline{H(Z,\nu)}.$$

Finally, we define the {\em local density} $d_0(X)$ of the space $X$
by $$d_0(X) = \min\{d(U) : U \in \tau^*(X)\}.$$ Clearly, we have
$$d_0(X) = \min\{|A| : A \in \mathcal{F}(X)\} = \min \{\Delta(D) : D \in \mathcal{D}(X)\}.$$
\end{definition}

The following result is obvious but very useful.

\begin{lemma}\label{lm:la}

Let $X$ be a space and $\lambda$ a singular cardinal such that every
$D \in \mathcal{D}(X)$ is $(<\lambda)$-resolvable. Then $X$ is
$\lambda$-resolvable.

\end{lemma}

As an immediate consequence of lemma \ref{lm:la} and of the first
part of theorem \ref{tm:spread} we obtain that if $\lambda$ is
singular and $s(X) < \lambda \le d_0(X) $ then $X$ is
$\lambda$-resolvable. (Of course, here $s(X) < \lambda$ is
equivalent with $\widehat{s}(X) < \lambda$ or with $\ps(X) <
\lambda$.)

The following lemma shows that, under certain simple and natural
conditions, if a space $X$ is not $\mu$-resolvable for some
 cardinal $\mu$ then some open set $V \in \tau^*(X)$
satisfies a condition just slightly weaker than $\mu \le d_0(V)$.

\begin{lemma}\label{lm:pm}
Let $X$ and ${\mu}$ be such that $c(X) < {\mu}\le
{\operatorname{\Delta}}(X)$. Then either $X$ is ${\mu}$-resolvable
or
\begin{enumerate}[(I)]\addtocounter{enumi}{1}
\item[$(*)$]
there is $V\in {\tau}^*(X)$ such that for each ${\kappa}<{\mu}$
there is $T \in \br V;<{\mu};$ with $d_0(V \backslash T) > \kappa$.
\end{enumerate}
If $\mu$ is regular then $V \in \tau^*(X)$ and $T \in [V]^{<\mu}$
may even be chosen so that $d_0(V \backslash T) \ge \mu$.
\end{lemma}

\begin{proof}[Proof of lemma \ref{lm:pm}]\prlabel{lm:pm}

Let us first consider the case when $\mu$ is regular and assume that
for all $V \in \tau^*(X)$ and $T \in [V]^{<\mu}$ we have $d_0(V
\backslash T) < \mu$. We define pairwise disjoint dense sets
$D_\alpha \in \mathcal{D}(X) \cap [X]^{<\mu}$ for $\alpha < \mu$ by
transfinite recursion as follows.

Assume that $\{D_\beta : \beta \in \alpha\} \subs \mathcal{D}(X)
\cap [X]^{<\mu}$ have already been defined and set $T = \cup
\{D_\beta : \beta \in \alpha\}$, then $|T| < \mu$ as $\mu$ is
regular. Let $\mathcal{W}$ be a maximal disjoint collection of open
sets $W \in \tau^*(X)$ such that $d(W \backslash T) < \mu$. By our
assumption, then $\cup \mathcal{W}$ is dense in $X$ and hence so is
$\cup \{W \backslash T : W \in \mathcal{W}\}$. So if for each $W \in
\mathcal{W}$ we fix $D_W \in \mathcal{D}(W \backslash T)$ with
$|D_W| < \mu$ then $D_\alpha = \cup \{D_W : W \in \mathcal{W}\}$ is
dense in $X$ as well and clearly $|D_\alpha| < \mu$. The family
$\{D_\alpha : \alpha < \mu\}$ witnesses that $X$ is
$\mu$-resolvable.

So let us assume now that $\mu$ is singular and fix   a strictly
increasing sequence $\<{\mu}_{\alpha}:{\alpha}<\cf({\mu})\>$ of
regular cardinals converging to  ${\mu}$ with $c(X) \cdot \cf(\mu) <
\mu_0$.

We then define a $\cf(\mu) \times \mu$ type matrix
$\{\aa{\alpha}{\xi}:{\alpha}<\cf({\mu}),{\xi}<{\mu}\}$ of pairwise
disjoint subsets of $X$, column by column in $\cf(\mu)$ steps, as
follows:
\begin{gather}\notag
X_{\alpha}=X\setm\bigcup\{\aa {\beta}{\xi}\ :\ {\beta}<{\alpha},\ 
{\xi}<{\mu}\},
\\\notag \aa{\alpha}{\xi}=\hh(X_{\alpha} \setm \cup \{
\aa{\alpha}{\zeta}\ :\ \zeta < \xi\},\ {\mu}_{\alpha}).
\end{gather}
Observe that we have $|\aa{\alpha}{\xi}| \le \mu_\alpha$, moreover
\begin{equation}\tag{$\dag$}\label{eq:dec}
\cl{\aa {\alpha}{\xi}}\supseteq\cl{\aa{\alpha}{\eta}} \text{
whenever   ${\alpha}<\cf({\mu})$ and ${\xi}\le {\eta}<{\mu}$}.
\end{equation}

Let us put $\aaa{\xi} = \bigcup\{\aa {\alpha}{\xi} : \alpha <
\cf(\mu)\}$ for ${\xi}<{\mu}$. The sets $\aaa {\xi}$ are pairwise
disjoint, so if they are all dense in $X$ then $X$ is
$\mu$-resolvable. Thus we can assume that at least one of them is
not dense in $X$, hence there is a nonempty open set $V\subs X$ and
an ordinal ${\xi}^\star<{\mu}$ such that $V\cap
\aaa{\xi^\star}=\empt$. Then we also have
\begin{equation}\tag{$\ddag$}\label{eq:disj}
\text{ $V\cap \aaa{\eta}=\empt$ for each ${\eta}\ge {\xi}^\star$ }
\end{equation}
 because
of (\ref{eq:dec}).

For ${\kappa}<{\mu}$ pick ${\beta}<\cf({\mu})$ with ${\kappa}\le
{\mu}_{\beta}$ and put
\begin{displaymath}
T = \bigcup\{\aa {\alpha}{\xi}:{\alpha}\le {\beta}, {\xi}<{\xi}^*\}.
\end{displaymath}
Then $|T| \le \mu_\beta \cdot |\xi^*| < \mu$ and it is immediate
from our definitions that then we have
$$d_0(V \backslash T) > \mu_\beta \ge \kappa.$$
\end{proof}

Before giving our next result we introduce a refined version of the
family of fat sets $\mathcal{H}(Y,\nu)$ defined above and of the
associated operator $H(Y,\nu)$. If a cardinal $\varrho < \nu$ is
also given, then we let $$\mathcal{H}(Y, \varrho , \nu) = \{A \in
 \mathcal{H}(Y, \nu) : \Delta(A) \ge \varrho\}.$$ Again, if $c(X) \le
 \nu$ and $\mathcal{H}(Y, \varrho\ ,\ \nu)$ is non-empty then $\mathcal{H}(Y, \varrho ,
 \nu)$ has a member $H(Y, \rho, \nu)$ of maximal closure. (If
$\mathcal{H}(Y,\varrho, \nu)$ is empty then we set $H(Y, \varrho,
\nu) = \emptyset$.)

\begin{lemma}\label{lm:red}
Assume that $X$ is a topological space and ${\mu}$ is a singular
cardinal with  $c(X) < {\mu}\le {\operatorname{\Delta}}(X)$,
moreover $X$ satisfies condition $(*)$ from lemma \ref{lm:pm}, i. e.
for every $\kappa < \mu$ there is a set $T \in [X]^{<\mu}$ such that
$d_0(X \backslash T) > \kappa$. Then we have  either (i) or (ii)
below.
\begin{enumerate}[(i)]
\item There is a disjoint family $\{D_\alpha : \alpha < \cf(\mu)\} \subs
\mathcal{F}(X) \cap [X]^{< \mu}$ such that $\Delta(D_\alpha)$
converges to $\mu$, moreover $$ \cup \{D_\gamma : \gamma \ge
\alpha\} \in \mathcal{D}(X)$$ for all $\alpha < \cf(\mu)$.
\item There are an open set $W \in \tau^*(X)$ and a set $T \in
[X]^{<\mu}$ with $d_0(W \backslash T) \ge \mu$.
\end{enumerate}
 \end{lemma}

\begin{proof}[Proof of \ref{lm:red}]\prlabel{lm:red}
Fix the same strictly increasing sequence
$\<{\mu}_{\alpha}:{\alpha}<\cf({\mu})\>$ of regular cardinals
converging to  ${\mu}$ with $c(X) \cdot \cf(\mu) < \mu_0$ as in the
above proof. Note that then for each $\alpha < \cf(\mu)$ we have
$$\mu_\alpha^- = \sup\{\mu_\beta : \beta < \alpha\} < \mu_\alpha.$$
Then by a straight-forward transfinite recursion on $\alpha <
\cf(\mu)$ we define disjoint sets $D_\alpha \in [X]^{<\mu}$  as
follows.

If $D_\beta$ has been defined for each $\beta < \alpha$ then set
$$D_\alpha = H(X \backslash \cup \{D_\beta : \beta < \alpha\},\ \mu_\alpha^-,\ \mu_\alpha).$$
(Note that $D_\alpha$ may be empty but it is a member of
$\mathcal{F}(X)$ if it is not.) Next, for each $\alpha < \cf(\mu)$
we let
$$E_\alpha = \cup \{D_\gamma : \gamma \ge \alpha\}.$$

Assume first that $E_\alpha \in \mathcal{D}(X)$ for all $\alpha <
\cf(\mu)$. In particular, then $D_\alpha \ne \emptyset$ for
cofinally many $\alpha < \cf(\mu)$, hence by re-indexing we may
actually assume that $D_\alpha \ne \emptyset$ for all $\alpha <
\cf(\mu).$ Now, $\Delta(D_\alpha) > \mu_\alpha^-$ immediately
implies that $\Delta(D_\alpha)$ converges to $\mu$, hence $(i)$ is
satisfied.

Next, assume that some $E_\alpha$ is not dense, hence there is a $W
\in \tau^*(X)$ with $W \cap E_\alpha = \emptyset$. Since $X$
satisfies $(*)$ there is a set $S \in [X]^{< \mu}$ such that $d_0(X
\backslash S) > \mu_\alpha.$ Let us set $$T = \cup \{D_\beta : \beta
< \alpha\} \cup S,$$ then $|T| < \mu$ as well, moreover we claim
that $d_0(W \backslash T) = \kappa \ge \mu$.

Assume, indirectly, that $U \in \tau^*(W)$ and $d(U \backslash T) =
\kappa < \mu.$ Since $U \backslash T \subs X \backslash S$ we have
$\kappa > \mu_\alpha,$ hence if $\delta < \cf(\mu)$ is chosen so
that $$\mu_\delta^- \le \kappa < \mu_\delta $$ then $\alpha <
\delta.$ Let $A$ be any dense subset of $U \backslash T$ of size
$\kappa$, then clearly $\Delta(A) = \kappa$ as well, moreover $A
\subs X \backslash \cup\{D_\beta : \beta < \delta\}$ holds because
$W \cap E_\alpha = \emptyset$. But then, by our definition, we have
$$A \in \mathcal{H}(X \backslash \cup\{D_\beta : \beta < \delta\}, \mu_\delta^-, \mu_\delta),$$
hence $A \subs \overline{D_\delta}$, contradicting that $W \cap
\overline{D_\delta} = \emptyset.$
\end{proof}

We now give one more easy result that, for a limit cardinal
$\lambda$, may be used to conclude $\lambda$-resolvability.

\begin{lemma}\label{lm:di}
Let $X$ be a space and $\lambda$ a limit cardinal and assume that
$\{D_\alpha : \alpha < \cf(\lambda)\}$ are {\em disjoint} subsets of
$X$ such that $$\cup \{D_\alpha : \beta \le \alpha < \cf(\lambda)\}
\in \mathcal{D}(X)$$ for every $\beta < \cf(\lambda).$ Assume also
that $D_\alpha$ is $\kappa_\alpha$-resolvable for each $\alpha <
\cf(\lambda)$ and the sequence $\langle \kappa_\alpha : \alpha <
\cf(\lambda)\rangle$ converges to $\lambda$. Then $X$ is
$\lambda$-resolvable.
\end{lemma}
\begin{proof}[Proof of \ref{lm:di}]
For each $\alpha < \cf(\lambda)$ fix a disjoint family
$$\{E_\xi^\alpha : \xi < \kappa_\alpha\} \subs \mathcal{D}(D_\alpha),$$
then for any $\xi < \lambda$ set $$E_\xi = \cup\{E_\xi^\alpha : \xi
< \kappa_\alpha\}.$$ Since the $\kappa_\alpha$ converge to
$\lambda$, for any fixed $\xi < \lambda$ we eventually have $\xi <
\kappa_\alpha$ and so $E_\xi$ is dense in $X$. Consequently the
disjoint family $\{E_\xi : \xi < \lambda\}$ witnesses that $X$ is
$\lambda$-resolvable.
\end{proof}

>From the above results and the first part of theorem \ref{tm:spread}
we may now easily obtain the "missing" second part. Indeed, assume
that $\lambda$ is singular and $s(X) < \Delta(X) = \lambda.$
Reasoning inductively, we may assume that if $s(Y) < \Delta(Y) <
\lambda$ then $Y$ is maximally, that is $\Delta(Y)$-resolvable.

Now, by lemma \ref{lm:elkin}, to prove that $X$ is
$\lambda$-resolvable it suffices to show that some subspace of $X$
is. Since $c(X) \le s(X)$, from lemmas \ref{lm:pm} and \ref{lm:red}
it follows that, if $X$ itself is not $\lambda$-resolvable, then
either there are a $W \in \tau^*(X)$ and a $T \in [W]^{<\lambda}$
such that $d_0(W \backslash T) \ge \lambda$ or there is a $V \in
\tau^*(X)$ with disjoint sets $\{D_\alpha : \alpha < \cf(\lambda)\}
\subs \mathcal{F}(V)$ such that $\Delta(D_\alpha)$ converges to
$\lambda$ and
$$ \cup \{D_\gamma : \alpha \le \gamma < \cf(\lambda) \} \in \mathcal{D}(X)$$ for
all $\alpha < \cf(\lambda)$. But we have seen that in the first case
$W \backslash T$ (and hence $W$), while in the second $V$ is
$\lambda$-resolvable.

We are now ready to present our result that, under certain
conditions, enables us to deduce $\lambda$-resolvability from
$(<\lambda)$-resolvability for a singular cardinal $\lambda$. We
first recall that $\celh(X)$ is defined as the smallest
(uncountable) cardinal such that $X$ has no disjoint family of open
sets of that size. As was shown in \cite{ET} (see also
\cite[4.1]{J}), $\celh(X)$ is always a regular cardinal. We also
note that if $\lambda$ is a limit cardinal then every
$(<\lambda)$-resolvable space $S$ has dispersion character
$\Delta(S) \ge \lambda$.

\begin{theorem}\label{tm:res}
Assume that $X$ is a topological space, $\lambda$ is a singular
cardinal, and  $\celh(X)\le \cf({\lambda}) <{\lambda}\le
{\operatorname{\Delta}}(X)$. If every dense subspace $S \subs X$
satisfying $\operatorname{\Delta}(S) \ge \lambda $ is
$(<\lambda)$-resolvable then $X$ is actually $\lambda$-resolvable.
\end{theorem}

\begin{proof}[Proof of \ref{tm:res}]\prlabel{tm:res}
Let us start by pointing out that if $A$ is fat in $X$ then $S = A
\cup (X \backslash \overline{A}) \in \mathcal{D}(X)$, moreover
$\Delta(A) \ge \lambda$ implies $\Delta(S) \ge \lambda$. So, every
fat set  $A \in \mathcal{F}(X)$ that satisfies $\Delta(A) \ge
\lambda$ is $(<\lambda)$-resolvable. It immediately follows from
this that the conditions on our space $X$ are inherited by all
non-empty open subspaces, hence by lemma \ref{lm:elkin} it is again
sufficient to prove that $X$ has some $\lambda$-resolvable subspace.

Now, if some $A \in \mathcal{F}(X)$ satisfies $d_0(A) \ge \lambda$
then $\Delta(B) \ge \lambda$ holds for every $B \in \mathcal{D}(A)$,
hence all dense subsets of $A$ are $(<\lambda)$-resolvable. But
then, by lemma \ref{lm:la}, $A$ is $\lambda$-resolvable.

Therefore, from here on we may assume that $d_0(A) < \lambda$ for
all $A \in \mathcal{F}(X).$ Actually, we claim that then even $d(A)
< \lambda $ holds whenever $A \in \mathcal{F}(X)$. Indeed, if $A \in
\mathcal{D}(U)$ for some $U \in \tau^*(X)$ then let $\mathcal{W}$ be
a maximal disjoint family of open sets $W \subs U$ such that $d(A
\cap W) < \lambda$. Then $\celh(X) \le \cf(\lambda) = \kappa$
implies $|\mathcal{W}| < \kappa$, moreover $\cup \mathcal{W}$ is
clearly dense in $U$ by our assumption. But then $\cup \mathcal{W}
\cap A$ is dense in $A$ and so $$d(A) \le d(\cup \mathcal{W} \cap A)
= \sum \{d(W \cap A) : W \in \mathcal{W}\} < \lambda.$$ (We note
that this is the only part of the proof where $\celh(X) \le
\cf(\lambda)$ is used rather than the weaker assumption $c(X) <
\lambda$.)

By lemma \ref{lm:pm}, if $X$ itself is not $\lambda$-resolvable then
there is a $V \in \tau^*(X)$ that satisfies condition $(*)$. We
shall show that then $V$ is $\lambda$-resolvable.

To see this, first fix a strictly increasing sequence $\langle
\lambda_\alpha : \alpha < \kappa \rangle$ of cardinals converging to
$\lambda$ and then, using $(*)$, fix for each $\alpha < \kappa$ a
set $T_\alpha \in [V]^{<\lambda}$ with $d_0(V \backslash T_\alpha) >
\lambda_\alpha$. Having done this, we define disjoint sets $D_\alpha
\in \mathcal{D}(V) \cap [V]^{<\lambda}$ by transfinite induction on
$\alpha < \kappa$ as follows.

Assume that $\alpha < \kappa$ and $D_\beta \in \mathcal{D}(V) \cap
[V]^{<\lambda}$ has been defined for each $\beta < \alpha$. Set
$$Z_\alpha = X \setminus (\cup \{D_\beta : \beta < \alpha\} \cup T_\alpha),$$
then $Z_\alpha$ is dense in $V$ because $\Delta(X) \ge \lambda.$ But
then $d(Z_\alpha) < \lambda$, hence we may pick $D_\alpha \in
\mathcal{D}(Z_\alpha) \subs \mathcal{D}(V)$ with $|D_\alpha| <
\lambda.$ Note that as $D_\alpha \subs V \setminus T_\alpha$ we also
have $\Delta(D_\alpha) > \lambda_\alpha.$

Now consider any partition $\{J_\xi : \xi < \kappa\}$ of $\kappa$
into $\kappa$ many sets of size $\kappa$ and for each $\xi < \kappa$
put
$$E_\xi = \cup \{D_\alpha : \alpha \in J_\xi\}.$$ Then each $E_\xi$
is dense in $V$ and clearly $\Delta(E_\xi) = \lambda$, hence it is
$(<\lambda)$-resolvable. But the $E_\xi$'s are pairwise disjoint,
hence obviously $V$ is $\lambda$-resolvable.
\end{proof}

We do not know if the assumption $\celh(X) \le \cf(\lambda)$ can be
relaxed to $c(X) < \lambda$ in theorem \ref{tm:res}, or even if it
can be dropped altogether.

\section{A simpler proof of Pavlov's theorem concerning extent}

The {\em extent} $\ee(X)$ of a space $X$ is defined as the supremum
of sizes of all closed discrete subspaces of $X$. (This is
Archangelski\v i's notation, in \cite{Pa} $\ext(X)$ and in \cite{J}
$\pp(X)$ is used to denote the same cardinal function.) Similarly as
in the previous section for the spread $\spre(X)$, we may define
$\operatorname{\widehat{e}}(X)$  as the smallest infinite (but not
necessarily uncountable) cardinal such that $X$ has no closed
discrete subset of that size. Note that a space $X$ is countably
compact iff $\operatorname{\widehat{e}}(X) = \omega.$ Clearly, one
has $\operatorname{\widehat{e}}(X) \le\pe(X)$ (the latter was
defined in the abstract).

In \cite{Pa} it was proved that $\Delta(X) > \pe(X)$ implies the
$\omega$-resolvability of $X$ for any $T_3$ space $X$. In this
section we shall present our proof of the slightly stronger result
in which only $\Delta(X) > \operatorname{\widehat{e}}(X)$ is used.  We believe that
this proof is significantly simpler than the one given in \cite{Pa},
although it follows the same steps.

We start with giving our simplified proof of the following result of
Pavlov concerning spaces that are finite unions of left separated
subspaces.

\begin{theorem}{\sc(Pavlov){\cite[Lemma 3.1]{Pa}}}\label{pavlov:left}
Assume that $\ls(X)<{\omega}$ and ${\kappa}\le |X|$ is an
uncountable regular cardinal. Then there is a strictly increasing
and continuous sequence $\langle F_{\alpha}:{\alpha}<{\kappa}
\rangle$ of closed subsets of $X$ with $|F_{\alpha}|<{\kappa}$ for
all $\alpha < \kappa$.
\end{theorem}

\begin{proof}
We prove the theorem by induction on $\ls(X)$. So assume that it is
true for $\ls(X) = k$ and consider  $X = \bigcup_{0 \le i \le k}
L_i$ where the $L_i$ are disjoint and left separated, moreover
$\omega < \kappa \le |X|.$ We may clearly assume that the left
separating order type of each $L_i$ is $\le \kappa$.

Assume that $S$ is an initial segment of some $L_i$ with $\tip(S) <
\kappa$ and $|\overline{S}| \ge \kappa$ (closures are always taken
in $X$). Since $\overline{S} \cap L_i = S$ we may apply the
inductive hypothesis to $\overline{S} \backslash S$ and find an
increasing and continuous $\kappa$-sequence $\langle
F_{\alpha}:{\alpha}<{\kappa} \rangle$ of its closed subsets of size
$< \kappa$. But then the traces $\overline{F_\alpha} \cap S$ will
stabilize and $|\overline{F_\alpha}| \le |F_\alpha| + |S| < \kappa$,
hence a suitable final segment of $\langle
\overline{F_{\alpha}}:{\alpha}<{\kappa} \rangle$ is as required.
Almost the same argument shows that the inductive step can also be
completed if $|L_i| < \kappa$ for some $i$. So we may assume that
$\tip{L_i} = \kappa$ for each $i$ and that $|\overline{A}| < \kappa$
whenever $A \in [X]^{<\kappa}$.

Let $y_\alpha$ denote the $\alpha$th member of $L_0$ and use the
inductive assumption to find an increasing and continuous
$\kappa$-sequence $\langle F_{\alpha}:{\alpha}<{\kappa} \rangle$ of
closed subsets of $ \bigcup_{1 \le i \le k} L_i$ of size $< \kappa$,
and then consider the set $$I = \{\alpha < \kappa : y_\alpha \in
\overline{F_\alpha}\}.$$

Assume first that $|I| < \kappa$ and hence $\sigma = \sup I <
\kappa$. We claim that then the set
\begin{displaymath}
J=\{{\beta} > \sigma :{\overline{F_{\beta}}} \neq
\cup_{{\gamma}<{\beta}} {\overline{F_{\gamma}}} \}
\end{displaymath}
is non-stationary in $\kappa.$ Indeed, for each $\beta \in J$ there
must be some $g(\beta) < \kappa$ with
$y_{g({\beta})}\in{\overline{F_{\beta}}}
\setm\cup_{{\gamma}<{\beta}} {\overline{F_{\gamma}}}$. Since
$g(\beta) \ge \beta > \sigma$ would imply $g(\beta) \notin I$ and
hence $$y_{g(\beta)} \notin \overline{F_{g(\beta)}} \supset
\overline{F_\beta},$$ we must have $g(\beta) < \beta$. But the
regressive function $g$ is clearly one-to-one on $J$, hence by
Fodor's (or Neumer's) pressing down theorem $J$ is non-stationary.
So there is a club set $C$ in $\kappa$ with $C \cap J = \emptyset$,
and then the sequence $\<{\overline{F_{\alpha}}}:{\alpha}\in C
\backslash \sigma \>$ clearly satisfies our requirements.

So we may assume that $|I| = \kappa.$ For each ${\alpha}<{\kappa}$
let us put $H_{\alpha}= \overline{\{y_{\gamma}:{\gamma}\in I\cap
{\alpha}\}}$. Note that we have $H_\alpha \subs \overline{F_\alpha}$
by the definition of $I$. Next, consider the set
\begin{displaymath}
J=\{{\alpha} < \kappa : \alpha \mbox{ is limit and } H_{\alpha} \neq
\cup_{{\gamma}<{\alpha}} H_{\gamma} \}.
\end{displaymath}
We claim that this set $J$ is again non-stationary. Indeed, for
every ${\alpha}\in J$ we may pick a "witness" $z_{\alpha}\in
H_{\alpha} \setm \cup_{{\gamma}<{\alpha}} H_{\gamma}$. Now, if
$z_\alpha \in L_0$ then $z_\alpha = y_{g(\alpha)}$ for some
$g(\alpha) < \alpha$ because $L_0$ is left separated. If, on the
other hand, $z_\alpha \notin L_0$ then $z_\alpha \in H_\alpha \subs
\overline{F_\alpha}$ implies $z_\alpha \in F_\alpha$ because
$F_\alpha$ is closed in $X \backslash L_0$. But the sequence
$\<F_{\alpha}:{\alpha}\in {\kappa}\>$ is continuous, hence in this
case we can choose an ordinal $g({\alpha})<{\alpha}$ such that
$z_{\alpha}\in F_{g({\alpha})}$.

In other words, this means that if $g(\alpha) = \beta$ then
$z_\alpha \in \{y_\beta\} \cup F_\beta.$ Now, the sequence $\langle
z_\alpha : \alpha \in J\rangle$ is obviously one-to-one, hence for
each ${\beta} < \kappa$ we have $|g^{-1}\{{\beta}\}|\le |F_{\beta}|
+ 1 <{\kappa}$, consequently, again by Fodor, $J$ is not stationary.
So there is a club $C\subs {\kappa}\setm J$ and then
$\<H_{\alpha}:{\alpha}\in C\>$ is increasing and continuous, however
maybe it is not strictly increasing. But $|I| = \kappa$ clearly
implies that the union of the $H_\alpha$'s is of size $\kappa$ and
so an appropriate subsequence of $\<H_{\alpha}:{\alpha}\in C\>$ will
be both continuous and strictly increasing.
\end{proof}

Before proceeding further, we need a simple definition.
\begin{definition}
Let $X$ be a space and $\mu$ an infinite cardinal number. We say
that $x \in X$ is a $T_\mu$ point of $X$ if for every set $A \in
[X]^{<\mu}$ there is some $B \in [X \backslash A]^{<\mu}$ such that
$x \in \overline{B}$. We shall use $T_\mu(X)$ to denote the set of
all $T_\mu$ points of $X$.
\end{definition}
For the reader familiar with Pavlov's paper \cite{Pa} we note that
his tr$_{\nu^+,\,\, \nu }(X)$ is identical with our $T_{\nu^+}(X)$.
Note also that if $Y \subs X$ then trivially any $T_\mu$ point in
$Y$ is a $T_\mu$ point in $X$, that is, we have $T_\mu(Y) \subs
T_\mu(X).$ Finally, if $\mu$ is regular then the set $T_\mu(X)$ is
clearly $(<\mu)$-closed in $X$, i. e. for every set $A \in
[T_\mu(X)]^{<\mu}$ we have $\overline{A} \subs T_\mu(X). $

\begin{lemma}\label{lm:dis}
Assume that the space $X$ may be written as the union of a strictly
increasing continuous chain $\langle F_{\alpha}:{\alpha}<{\kappa}
\rangle$ of closed subsets of $X$ with $|F_{\alpha}|<{\kappa}$ for
all $\alpha < \kappa$, where $\kappa$ is an uncountable regular
cardinal. Then $T_\kappa(X) = \emptyset$ implies that there exists a
set $D \subs X$ with $|D| = \kappa$ such that every subset $Y \in
[D]^{<\kappa}$ is closed discrete in $X$. In particular, we have
$\operatorname{\widehat{e}}(X) \ge \kappa.$
\end{lemma}

\begin{proof}
The assumption $T_\kappa(X) = \emptyset$ implies that for every
point $x \in X$ we may fix a set $A_x \in [X]^{<\kappa}$ such that
$x \notin \overline{B}$ whenever $B \in [X \backslash
A_x]^{<\kappa}$. By the regularity of $\kappa$, the set
$$C = \{\alpha < \kappa : \forall x \in F_\alpha  (A_x \subs
F_\alpha) \}
$$ is club in $\kappa$. For each $\alpha \in C$ let us pick a point
$x_\alpha \in F_{\alpha + 1} \backslash F_\alpha$ and then set $D =
\{x_\alpha : \alpha \in C\}.$

To see that this $D$ is as required, it remains to show that all
"small" subsets of $D$ are closed discrete. This in turn will follow
if we show that all proper initial segments of $D$ are. So let
$\gamma < \kappa$ and consider the set $S = \{x_\alpha : \alpha \in
C \cap \gamma\}$. For every point $y \in X$ there is a $\beta <
\kappa$ such that $y \in F_{\beta+1} \backslash F_\beta$. Let
$\delta$ be the largest element of $C$ with $\delta \le \beta$ and
$\varepsilon$ the smallest element of $C$ above $\beta$, hence we
have $\delta \le \beta < \varepsilon.$

Then, on one hand, $\{x_\alpha : \alpha < \delta\} \subs F_\delta
\subs F_\beta$, while on the other hand $A_y \subs F_\varepsilon$
and $\{x_\alpha : \varepsilon \le \alpha < \gamma\} \subs X
\backslash F_\epsilon$, which together imply that $y$ has a
neighbourhood $U$ such that $U \cap S \subs \{x_\delta\}$.
\end{proof}

We need one more result making use of the operator $T_\mu.$

\begin{lemma}\label{lm:tm}
If a space $X$ satisfies $T_\mu(X) = X$ for a regular cardinal $\mu$
then $X$ is $\mu$-resolvable.
\end{lemma}

\begin{proof}
Clearly, $T_\mu(X) = X$ implies $T_\mu(U) = U$ for all open subsets
$U \subs X$, hence by lemma \ref{lm:elkin} it suffices to show that
$X$ includes a $\mu$-resolvable subspace $Y$.

Since every point of $X$ is a $T_\mu$ point, for any set $A \in
[X]^{<\mu}$ we may fix a {\em disjoint} family $\bcal(A) \subs [X
\backslash A]^{<\mu}$ with $|\bcal(A)| = |A| < \mu$ such that  $$A
\subs \cup \{\overline{B} : B \in \bcal(A)\}.$$

We now define sets $A_\alpha$ in $[X]^{<\mu}$ by induction on
$\alpha < \mu$ as follows. Let $x \in X$ be any point and start with
$A_0 = \{x\}$. Assume next that $0 < \alpha < \mu$ and the sets
$A_\beta \in [X]^{<\mu}$ have been defined for all $\beta < \alpha.$
Then we set
$$B_\alpha = \bigcup \bcal\big(\cup \{A_\beta : \beta < \alpha\}\big)
\mbox{ and } A_\alpha = \cup \{A_\beta : \beta < \alpha\} \cup
B_\alpha.$$ After the induction is completed we let $$Y = \cup
\{A_\alpha : \alpha < \mu\}.$$

It is clear from the construction that the $B_\alpha$'s are pairwise
disjoint, moreover for every set $s \in [\mu]^\mu$ the union
$\cup_{\alpha \in s}B_\alpha$ is dense in $Y$. But then $Y$ is
obviously $\mu$-resolvable.

\end{proof}

We are now ready to state and prove our promised result.

\begin{theorem}\label{tm:sim}
Assume that the regular closed subsets of the space $X$ form a
$\pi$-network in $X$ and $T_\mu(X)$ is dense in $X$ for some regular
cardinal $\mu> \operatorname{\widehat{e}}(X)$.  Then $X$ is
$\omega$-resolvable. In particular, any $T_3$ space $X$ satisfying
$\Delta(X) > \operatorname{\widehat{e}}(X)$ is $\omega$-resolvable.
\end{theorem}

\begin{proof}
Assume, indirectly, that $X$ is $\omega$-irresolvable. By lemmas
\ref{lm:elkin} and \ref{lm:ls} then there is a regular closed subset
$K$ of $X$ that is both {\em hereditarily} $\omega$-irresolvable and
satisfies $\ls(K) < \omega.$

Let us now define the sequence of sets $\{K_n : n <\omega\}$ by the
following recursion: $K_0 = K$ and $K_{n+1} = T_\mu(K_n)$. Since
$T_\mu(Y)$ is $(<\mu)$-closed in $Y$ for any space $Y$, we may
conclude by a simple induction that $K_i$ is $(<\mu)$-closed in $K$
and hence $\operatorname{\widehat{e}}(K_i) \le \operatorname{\widehat{e}}(K) < \mu$ for all $i <
\omega$.

We next claim that, for each $n < \omega$, $K_{n+1} = T_\mu(K_n)$ is
dense in $K_n$ and hence in $K$. For $n = 0$ this follows
immediately from our assumption that $T_\mu(X) \in \dcal(X).$

Clearly, any neighborhood of a $T_\mu$ point in any space must have
size at least $\mu$. Hence if our claim holds up to (and including)
$n$ then we also have $\Delta(K_n) \ge \mu$ and since $K_n \in
\dcal(K)$ the regular closed subsets of $K_n$ form a $\pi$-network
in $K_n$. (The latter holds because the regular closed subsets of a
dense subspace are exactly the traces of the regular closed sets in
the original space.)

 Now, let $U$ be any non-empty open subset of $K_n$.
We show first that then $|U \cap K_{n+1}| \ge \mu$, hence
$\Delta(K_{n+1}) \ge \mu$. (In other words, $K_{n+1}$ is not only
dense but even $\mu$-dense in $K_n$.) To see this, let $\emptyset
\ne F \subs U$ be regular closed in $K_n$, then $|F| \ge \mu$ and
$\ls(F) < \omega$ imply, in view of theorem \ref{pavlov:left}, the
existence of a strictly increasing continuous sequence $\langle
F_\alpha : \alpha < \mu
 \rangle$ of closed subsets of $F$ (and hence of $X$) with
$|F_\alpha| < \mu$. Then we may apply lemma \ref{lm:dis} to any
final segment of the sequence $\langle F_\alpha : \alpha < \mu
 \rangle$ to conclude that $F_\alpha \cap T_\mu(K_n) = F_\alpha \cap
K_{n+1} \ne \emptyset$ for cofinally many $\alpha < \mu$, hence $|U
\cap K_{n+1}| \ge |F \cap K_{n+1}| \ge \mu$.

But $\Delta(K_{n+1}) \ge \mu$ implies that for any non-empty regular
closed set $H$ in $K_{n+1}$ we have $|H| \ge \mu$, and so, using
again $\ls(H) < \omega$ and $\operatorname{\widehat{e}}(K_n) < \mu$,
we obtain from theorem \ref{pavlov:left} and lemma \ref{lm:dis} that
$T_\mu(H) $ is non-empty, i. e. $K_{n+2}$ is indeed dense in
$K_{n+1}$.

Now suppose that there is an $n < \omega$ such that $K_n \backslash
K_{n+1}$ is not dense in $K_n$. This would imply that for some $U
\in \tau^*(K_n)$ we have $U \subs K_{n+1}$ and hence $T_\mu(U) = U$.
But that would imply by lemma \ref{lm:tm} that $U$ is
$\mu$-resolvable, a contradiction. Therefore, we must have that $K_n
\backslash K_{n+1}$ is dense in $K_n$ and hence in $K$ for all $n <
\omega$. But then $K$ would be $\omega$-resolvable, which is again
absurd. This contradiction then completes the proof of the first
part of our theorem.

To see the second part note that, by lemma \ref{lm:ls} and by
considering regular closed subsets of $X$, it suffices to prove the
$\omega$-resolvability of $X$ under the additional condition $\ls(X)
< \omega$. But then $T_\mu(X) \in \dcal(X)$ follows immediately from
theorem \ref{pavlov:left} and lemma \ref{lm:dis} with the choice
$\mu = \operatorname{\widehat{e}}(X)^+$.
\end{proof}

Since for any crowded (i. e. dense-in-itself) countably compact
$T_3$ space $X$ one has $\Delta(X) \ge \mathfrak{c} \ge \omega_1$,
theorem \ref{tm:sim} immediately implies the following result of
Comfort and Garcia-Ferreira.

\begin{qtheorem}[Comfort,Garcia-Ferreira, {\cite[Theorem 6.9]{Co}}]
Every crowded and  countably compact $T_3$ space is
${\omega}$-resolvable.
  \end{qtheorem}

Note that the assumption of regularity in this theorem is essential
because of the following two results.

\begin{qtheorem}[Malykhin, {\cite[Example 14]{Ma}}]
There is a countably compact, irresolvable $T_2$ space.
\end{qtheorem}

\begin{qtheorem}[Pavlov, {\cite[Example 3.9]{Pa}}]
There is a countably compact, irresolvable Uryshon space.
\end{qtheorem}

Pytkeev has recently announced in \cite{Py2} that a crowded and
countably compact $T_3$ space is even ${\omega}_1$-resolvable. We
haven't seen his paper but would like to point out that this
stronger result is an immediate consequence of an old (and deep)
result of Tka\v cenko and of lemma \ref{lm:ls}.

Tka\v cenko proved in \cite {T1} that {\em if  $X$ is a countably
compact $T_3$ space with $\ls(X) \le \omega $ then $X$ is compact
and scattered.} In \cite{GJSz} it was shown that this statement
remains valid if $T_3$ is weakened to $T_2$, hence we get the
following result.
\begin{theorem}\label{tm:coco}
If $X$ is a crowded and countably compact $T_2$ space in which the
regular closed subsets form a $\pi$-network then $X$ is
$\omega_1$-resolvable.
\end{theorem}
\begin{proof}[Proof of theorem \ref{tm:coco}]\prlabel{tm:coco}
By the above result from \cite{GJSz}, every non-empty regular closed
subset $F \subs X$ must satisfy $\ls(F) \ge \omega_1$. But then $X$
is $\omega_1$ -resolvable by lemma \ref{lm:elkin}.
\end{proof}

Any crowded and countably compact $T_3$ space has dispersion
character $\ge \mathfrak{c}$. Hence the following interesting, and
apparently difficult, problem is left open by theorem \ref{tm:coco}.

\begin{problem}
Is a crowded and countably compact $T_3$ space
$\mathfrak{c}$-resolvable or even maximally resolvable?
\end{problem}


\begin{thebibliography}{9}


\bibitem{Ba} {~K.~P.~S. Bhaskara Rao,}
{\em On $\aleph$-resolvability}, preprint.



\bibitem{Co}
 W.W. Comfort and S. Garcia-Ferreira,
{\em Resolvability: A selective survey and some new results,} Topology
and its Applications 74 (1996) 149--167.

\bibitem{El}
{A. G. El'kin, } {\em Resolvable spaces which are not maximally
resolvable.} Vestnik Moskov. Univ. Ser. I Mat. Meh. 24 1969 no. 4
66--70.

\bibitem{ET}
P. Erd\H os and A. Tarski, {\em On families of mutually exclusive
sets,} Ann. of Math. 44 (1943) 315--329.

\bibitem{GJSz} {J.  Gerlits, I. Juh\'asz, and Z. Szentmikl\'ossy}
{\em Two improvements on Tka\v cenko's addition theorem}, CMUC 46
(2005), 705-710.


\bibitem{HJ} A. Hajnal and  I. Juh\'asz,
{\em Discrete subspaces of topological spaces. II.}
Indag. Math. 31 (1969) 18--30.


\bibitem{J}  I. Juh\'asz, {\em Cardinal functions -- ten years
later}, Math. Center Tract no. 123, Amsterdam, 1980


\bibitem{df} I. Juh\'asz, L. Soukup, and Z. Szentmikl\'ossy,
{\em $\dcal$-forced spaces: a new approach to resolvability}, Top.
Appl.

\bibitem{Ma}
 {V. I. Malykhin, }
{\em  Resolvability of $A$-, ${\rm CA}$- and ${\rm PCA}$-sets in compacta.}
{Topology Appl.} 80 (1997), no. 1-2, 161--167.

\bibitem{Pa} O. Pavlov, {\em On resolvability of topological spaces}
Topology and its Applications 126 (2002) 37-47.

\bibitem{Py} E. G. Pytkeev, {\em Maximally decomposable spaces}, Trudy
Math. Inst. Steklova, 154 (1983), 209-213.

\bibitem{Py2} E. G. Pytkeev, {\em
Resolvability of countably compact regular spaces},
Proc. Steklov Inst. Math. 2002, Alg., Top., Math. Anal., suppl. 2,
S152-S154.


\bibitem{T1} M. G. Tka\v cenko: {\it
  \begin{cyr}
O bikompaktah, predstavimych v vide obedineniya schetnogo chisla
levych podprostranstv
  \end{cyr}},
CMUC {\bf 20}(1979), 361-379 and 381-395.



\end{thebibliography}
\end{document}